\documentclass[11pt,leqno]{article}
\usepackage{amsmath, amsfonts, amssymb, amsthm, graphicx}

\numberwithin{equation}{section}

\newtheorem{theorem}{Theorem}
\newtheorem{lemma}{Lemma}
\newtheorem{remark}{Remark}

\def\R{\mathbb{R}}

\newenvironment{proof1}{
    \noindent {\it Proof }}{\hfill$\Box$
}

\begin{document}

\title{\textbf{Wave breaking in the short-pulse equation}}
\author{Yue Liu$^1$, Dmitry Pelinovsky$^2$,  and Anton Sakovich$^2$ \\
{\small $^{1}$ Department of Mathematics, University of Texas at Arlington, Arlington, TX, 76019, USA} \\
{\small $^{2}$ Department of Mathematics, McMaster
University, Hamilton, ON, L8S 4K1, Canada}}
\date{}
\maketitle

\begin{abstract}
Sufficient conditions for wave breaking are found for the
short-pulse equation describing wave packets of few cycles on the ultra-short
pulse scale. The
analysis relies on the method of characteristics and conserved
quantities of the short-pulse equation and holds both on an
infinite line and in a periodic domain. Numerical illustrations of
the finite-time wave breaking are given in a periodic domain.
\end{abstract}

\section{Introduction}

The short-pulse equation,
\begin{equation}
\label{short-pulse}
u_{tx} = u + \frac{1}{6} (u^3)_{xx}, \quad t > 0, \; x\in \mathbb{R},
\end{equation}
is a useful and simple approximation of nonlinear wave packets in
dispersive media in the limit of few cycles on the ultra-short
pulse scale \cite{CJSW05,ScWa}. This equation is a dispersive
generalization of the following advection equation
\begin{equation}
\label{simple-wave}
u_t = \frac{1}{2} u^2 u_x, \quad t > 0, \; x\in \mathbb{R}.
\end{equation}

According to the method of characteristics, the advection equation
(\ref{simple-wave}) exhibits wave breaking in a finite time for
any initial data $u(x,0) = u_0(x)$ on an infinite line if $u_0(x)$ is continuously differentiable and
there is a point $x_0 \in \R$ such that $u_0(x_0) u_0'(x_0) > 0$. This follows from the implicit
solution
$$
u(\xi,t) = u_0(\xi), \quad x(\xi,t) = \xi - \frac{1}{2} u_0^2(\xi) t, \quad t > 0, \quad \xi \in \R,
$$
for any given $u_0(x) \in C^1(\R)$. We say that the finite-time wave breaking occurs if
there exists a finite time $T \in (0,\infty)$ such that
\begin{equation}
\label{def-wave-breaking}
\lim_{t \uparrow T} \sup_{x \in \R} u(x,t) u_x(x,t) = \infty, \quad
\mbox{while} \quad \lim_{t \uparrow T} \sup_{x \in \R} |u(x,t)| <
\infty.
\end{equation}
For the simple advection equation (\ref{simple-wave}),
$$
T = \inf_{\xi \in \R} \left\{ \frac{1}{u_0(\xi) u_0'(\xi)} : \quad u_0(\xi) u_0'(\xi) > 0 \right\}.
$$
In view of this result, we address the question if the
dispersion term $\partial_x^{-1} u$ in the
short-pulse equation (\ref{short-pulse}) can stabilize global
dynamics of the advection equation (\ref{simple-wave}) at least
for small initial data. Local well-posedness of the short-pulse
equation on an infinite line was proven in \cite{ScWa}.

\begin{theorem}[Sch\"{a}fer \& Wayne, 2004]
Let $u_0 \in H^2(\mathbb{R})$. There exists a $T > 0$ such that
the short-pulse equation (\ref{short-pulse}) admits a unique
solution
$$
u(t) \in C([0,T),H^2(\mathbb{R})) \cap C^1([0,T),H^{1}(\mathbb{R}))
$$
satisfying $u(0) = u_0$. Furthermore, the solution $u(t)$ depends
continuously on $u_0$. \label{theorem-wayne}
\end{theorem}

To extend local solutions into a global solution, Pelinovsky \&
Sakovich \cite{PelSak} used the following conserved quantities of
the short-pulse equation:
\begin{eqnarray}
\label{conserved-quantities-1}
E_0 & := & \int_{\mathbb{R}} u^2 dx, \\ \label{conserved-quantities-2a}
E_1 & := & \int_{\mathbb{R}} \left( \sqrt{1 + u_x^2} -1 \right) dx =
\int_{\mathbb{R}} \frac{u_x^2}{1 + \sqrt{1 + u_x^2}} dx, \\
\label{conserved-quantities-2}
E_2 & := & \int_{\mathbb{R}} \sqrt{1 + u_x^2}
\left[ \partial_x \left( \frac{u_x}{\sqrt{1+u_x^2}} \right) \right]^2 dx =
\int_{\mathbb{R}} \frac{u_{xx}^2}{(1 + u_x^2)^{5/2}} dx.
\end{eqnarray}
If $u(t)$ is a local solution in Theorem \ref{theorem-wayne}, then
$E_0$, $E_1$, and $E_2$ are bounded and constant in time for all
$t \in [0,T)$. Global well-posedness of the short-pulse equation
on an infinite line was proven in \cite{PelSak} for small initial
data in $H^2$ satisfying
\begin{equation}
\label{well-posedness-bound} 2 E_1 + E_2 \leq \| u_0' \|^2_{L^2} +
\| u_0'' \|^2_{L^2} < 1.
\end{equation}
The condition (\ref{well-posedness-bound}) can be sharpen using
the scaling transformation for the short-pulse equation
(\ref{short-pulse}). Let $\alpha \in \mathbb{R}_+$ be an arbitrary
parameter. If $u(x,t)$ is a solution of (\ref{short-pulse}), then
$\tilde{u}(\tilde{x},\tilde{t})$ is also a solution of (\ref{short-pulse}) with
\begin{equation}
\label{scaling-invariance} \tilde{x} = \alpha x, \quad \tilde{t} = \alpha^{-1} t,
\quad \tilde{u}(\tilde{x},\tilde{t}) = \alpha u(x,t).
\end{equation}
The conserved quantities transform as follows:
\begin{eqnarray*}
\tilde{E}_1 & = & \int_{\mathbb{R}} \left( \sqrt{1 + \tilde{u}_{\tilde{x}}^2} -1
\right) d \tilde{x} =
\alpha \int_{\mathbb{R}} \left( \sqrt{1 + u_x^2} -1 \right) dx  = \alpha E_1, \\
\tilde{E}_2 & = & \int_{\mathbb{R}} \frac{\tilde{u}_{\tilde{x} \tilde{x}}^2}{(1+
\tilde{u}_{\tilde{x}}^2)^{5/2}} d \tilde{x} = \alpha^{-1} \int_{\mathbb{R}}
\frac{u_{xx}^2}{(1+u_x^2)^{5/2}} dx = \alpha^{-1} E_2.
\end{eqnarray*}
Finding the minimum of $2 \tilde{E}_1 + \tilde{E}_2 = 2 \alpha E_1 + \alpha^{-1} E_2$ in $\alpha$
gives a sharper sufficient condition on global well-posedness \cite{PelSak}.

\begin{theorem}[Pelinovsky \& Sakovich, 2009]
\label{theorem-wellposedness} Let $u_0 \in H^2(\mathbb{R})$ and $2
\sqrt{2 E_1 E_2} < 1$. Then the short-pulse equation
(\ref{short-pulse}) admits a unique global solution $u(t) \in
C(\mathbb{R}_+,H^2(\mathbb{R}))$ satisfying $u(0) = u_0$.
\end{theorem}

Theorem \ref{theorem-wellposedness} does not exclude wave breaking
in a finite time for large initial data and this paper gives a
proof that the wave breaking may occur in the short-pulse
equation (\ref{short-pulse}). Negating the sufficient condition
for global well-posedness in Theorem \ref{theorem-wellposedness},
a necessary condition for the wave breaking follows: the wave breaking
may occur in the short-pulse equation (\ref{short-pulse}) with
the initial data $u_0 \in H^2$ only if $2 \sqrt{2 E_1 E_2} \geq 1$.
We shall find a sufficient condition for the wave breaking in the short-pulse
equation (\ref{short-pulse}).

Unlike the previous work in \cite{PelSak}, we will not be using
conserved quantity $E_2$ but will rely on the conservation of $E_0$, $E_1$,
and the energy
\begin{eqnarray}
\label{conserved-quantities-3} E_{-1} := \int_{\mathbb{R}} \left[
\left( \partial_x^{-1} u \right)^2 - \frac{1}{12} u^4 \right] dx.
\end{eqnarray}
Here $\partial_x^{-1} u$ is defined from a local solution $u$ by
$$
\partial_x^{-1} u := \int_{-\infty}^x u(x',t) dx' = -\int_x^{\infty} u(x',t) dx' =
\frac{1}{2} \left( \int_{-\infty}^x - \int_x^{\infty} \right) u(x',t) dx',
$$
thanks to the zero-mass constraint $\int_{\mathbb{R}} u(x,t) dx = 0$ for all $t \in (0,T)$.
(Note that the initial data $u_0$ does not have generally to satisfy the zero-mass constraint
$\int_{\mathbb{R}} u_0(x) dx = 0$.)

Thanks to the Sobolev inequality, $\| u \|_{L^4} \leq C \| u \|_{H^1}$ for some $C > 0$,
the quantity $E_{-1}$ is bounded if $u \in H^2(\mathbb{R}) \cap \dot{H}^{-1}(\mathbb{R})$, where
$\dot{H}^{-1}$ is defined by its norm
$$
\| u \|_{\dot{H}^{-1}} := \| \partial_x^{-1} u \|_{L^2}.
$$
(Note that if $u \in H^2(\mathbb{R}) \cap
\dot{H}^{-1}(\mathbb{R})$, then $\int_{\mathbb{R}} u(x) dx = 0$
is satisfied.) Our main result on the wave breaking on an infinite
line is formulated as follows.

\begin{theorem}
\label{theorem-wave-breaking} Let $u_0 \in H^2(\mathbb{R}) \cap
\dot{H}^{-1}(\mathbb{R})$ and $T$ be the maximal existence time of
Theorem \ref{theorem-wayne}. Let
\begin{eqnarray*}
F_1 & := & \frac{1}{\sqrt{2}} \left( E_1^2 + \left( 8 E_0 E_1 + E_1^4 \right)^{1/2} \right)^{1/2}, \\
F_0 & := & \frac{1}{\sqrt{2}} \left( E_0 + E_{-1} +
\frac{1}{12} E_0 F_1^2 \right)^{1/2},
\end{eqnarray*}
and assume that there exists $x_0 \in \mathbb{R}$ such that $u_0(x_0) u_0'(x_0) > 0$ and
\begin{eqnarray*}
& \mbox{\rm either} & \quad |u_0'(x_0)| > \left( \frac{F_1^2}{4 F_0} \right)^{1/3}, \quad
|u_0(x_0)| |u_0'(x_0)|^2 > F_1 + \left( 2 F_0 |u_0'(x_0)|^3 -
\frac{1}{2} F_1^2 \right)^{1/2}, \\
& \mbox{\rm or} & \quad
|u_0'(x_0)| \leq \left( \frac{F_1^2}{4 F_0} \right)^{1/3}, \quad
|u_0(x_0)| |u_0'(x_0)|^2 > F_1.
\end{eqnarray*}
Then $T < \infty$, so that the solution $u(t) \in
C([0,T),H^2(\mathbb{R})\cap \dot{H}^{-1}(\mathbb{R}))$ of the
short-pulse equation (\ref{short-pulse}) blows up in the sense of
$$
\lim_{t \uparrow T} \sup_{x \in \mathbb{R}} u(x,t) u_x(x,t) =
\infty, \quad \mbox{while} \quad \lim_{t \uparrow T} \| u(\cdot,t)
\|_{L^{\infty}} \leq F_1.
$$
\end{theorem}

\begin{remark}
\label{remark-scaling-invariance}
The quantities $F_0$ and $F_1$ in Theorem \ref{theorem-wave-breaking} 
can be defined by 
$$
F_1 := \sup_{t \in [0,T)} \| u(\cdot,t) \|_{L^{\infty}}, \quad 
F_0 := \sup_{t \in [0,T)} \|  \partial_x^{-1} u(\cdot,t) \|_{L^{\infty}}.
$$
Note that the scaling transformation (\ref{scaling-invariance}) gives
$$
\| \tilde{u}(\cdot,\tilde{t}) \|_{L^{\infty}} = \alpha \| u(\cdot,t) \|_{L^{\infty}}, \quad
\| \partial_{\tilde{x}}^{-1} \tilde{u}(\cdot,\tilde{t}) \|_{L^{\infty}} = \alpha^2
\| \partial_{x}^{-1} u(\cdot,t) \|_{L^{\infty}},
$$
so that the sufficient condition
of Theorem \ref{theorem-wave-breaking} with new definitions of 
$F_1$ and $F_0$ is invariant in $\alpha$. We note, however, that, 
while the bound on $F_1$ in Theorem \ref{theorem-wave-breaking} 
scales correctly as $\tilde{F}_1 = \alpha F_1$,
the bound on $F_0$ is not correctly scaled in $\alpha$ because
$$
\tilde{E}_0 = \alpha^3 E_0, \quad \tilde{E}_{-1} = \alpha^5 E_{-1}.
$$
This is an artefact of using Sobolev embedding in Lemma \ref{lemma-bounds-F0} below.
\end{remark}

To prove Theorem \ref{theorem-wave-breaking}, we shall adopt the
method of characteristics and proceed with apriori differential
estimates. Our results remain valid in a periodic domain, where Theorem
\ref{theorem-wave-breaking} is replaced by Theorem \ref{theorem-breaking-periodic}
below. The technique of characteristics and apriori differential estimates has been
applied for wave breaking in other nonlinear wave equations, see
\cite{Constantin,C-E,Chae,Hu,LY1,LY2,Yin} for an incomplete list of references.

We emphasize that, unlike the previous work in \cite{PelSak}, we avoid using a
transformation between the short-pulse equation
(\ref{short-pulse}) and the integrable sine--Gordon equation in
characteristic coordinates. Our proof of the wave breaking for the
short-pulse equation (\ref{short-pulse}) does not suggest that
there exists a similar wave breaking for the sine--Gordon
equation, it is rather the breaking of the coordinate
transformation between the two equations. On a similar note, we do
not use the integrability properties of the short-pulse equation
(\ref{short-pulse}) such as the Lax pair, the inverse scattering
transform method, the bi-Hamiltonian formulation, and the
existence of exact soliton solutions.

The article is constructed as follows. The proof of Theorem \ref{theorem-wave-breaking} is
given in Section 2. Section 3 reports extension of Theorem \ref{theorem-wave-breaking} to
a periodic domain. Section 4 contains numerical evidences of the finite-time
wave breaking in a periodic domain.

\section{Wave breaking on an infinite line}

Let us rewrite the Cauchy problem for the short-pulse equation on
an infinite line in the form
\begin{equation}
\label{Cauchy}
\left\{ \begin{array}{l} u_t = \frac{1}{2} u^2 u_x + \partial_x^{-1} u, \quad x \in \mathbb{R}, \;\; t > 0 \\
u(x,0) = u_0(x), \quad \quad \;\; x \in \mathbb{R}, \end{array} \right.
\end{equation}
where $\partial_x^{-1} u := \int_{-\infty}^x u(x',t) dx'$. In what follows, we use
both notations $u(t)$ and $u(x,t)$ for the same solution of the Cauchy problem (\ref{Cauchy}).
Local existence of solutions with the conservation of $E_{-1}$ and $E_0$ is described by
the following result.

\begin{lemma}
\label{lemma-conservation}
Let $u_{0} \in H^s(\mathbb{R}) \cap \dot{H}^{-1}(\mathbb{R})$, $s \geq 2$.
There exist a maximal time $T = T(u_0) > 0$ and a unique solution $u(x,t)$
to the Cauchy problem (\ref{Cauchy}) such that
$$
u(t) \in C([0,T),H^s(\mathbb{R}) \cap \dot{H}^{-1}(\mathbb{R})) \cap
C^{1}([0,T),H^{s-1}(\mathbb{R}))
$$
satisfying $u(0) = u_0$. Moreover, the solution $u(t)$ depends continuously on the initial
data $u_0$ and the values of $E_{-1}$, $E_0$, and $E_1$
in (\ref{conserved-quantities-1}), (\ref{conserved-quantities-2a}), and (\ref{conserved-quantities-3})
are constant on $[0,T)$.
\end{lemma}

\begin{proof}
If $u_0 \in H^s(\mathbb{R}) \cap \dot{H}^{-1}(\mathbb{R})$, $s
\geq 2$, then $\partial_x^{-1} u_0 \in H^{s+1}(\mathbb{R})$, so
that $\int_{\mathbb{R}} u_0(x) dx = 0$. By the theorem of
Sch\"{a}fer \& Wayne \cite{ScWa}, there exists a solution
$$
u(t) \in C([0,T),H^s(\mathbb{R})) \cap C^{1}([0,T),H^{s-1}(\mathbb{R}))
$$
of the short-pulse equation
(\ref{short-pulse}), so that
$$
\partial_x^{-1} u(t) := u_t - \frac{1}{2} u^2 u_x \in
C((0,T),H^{s-1}(\mathbb{R})).
$$
Therefore, $u(t) \in C([0,T),H^s(\mathbb{R})\cap
\dot{H}^{-1}(\mathbb{R}))$ in view of continuity of $\| u
\|_{\dot{H}^{-1}}$ at $t = 0$. Because $f \in H^s(\mathbb{R})$, $s \geq 1$
implies $\lim_{|x| \to \infty} f(x) = 0$, the
zero-mass constraint holds in the form
\begin{equation}
\label{mass-conservation} \int_{\mathbb{R}} u(x,t) dx = 0, \quad
t \in [0,T).
\end{equation}
Let us define
$$
\partial_x^{-2} u(t) := \left(\partial_x^{-1} u \right)_t - \frac{1}{6} u^3.
$$
By the zero-mass constraint (\ref{mass-conservation}) and
uniqueness of the solution $u(t)$ for any $t \in [0,T)$, we obtain
$$
\lim_{|x| \to \infty} \partial_x^{-2} u(x,t) = 0, \quad t \in
[0,T).
$$
Using balance equations for the densities of $E_{-1}$, $E_0$, and $E_1$, we write
\begin{eqnarray*}
\left[ (\partial_x^{-1} u)^2 - \frac{1}{12} u^4 \right]_t & = & \left[
(\partial_x^{-2} u)^2 - \frac{1}{36} u^6 \right]_x, \\
\left( u^2 \right)_t & = & \left( (\partial_x^{-1} u)^2 + \frac{1}{4} u^4 \right)_x, \\
\left( \sqrt{1+u_x^2} - 1 \right)_t & = & \frac{1}{2}
\left( u^2 \sqrt{1 + u_x^2} \right)_x.
\end{eqnarray*}
Integrating the balance equation in $x \in \mathbb{R}$ for any $t \in
[0,T)$, we complete the proof that $E_{-1}$, $E_0$, and $E_1$ are bounded
and constant on $[0,T)$.
\end{proof}

\begin{remark}
The maximal existence time $T > 0$ in Lemma
\ref{lemma-conservation} is independent of $s \geq 2$ in the
following sense. If $u_0 \in H^s(\mathbb{R}) \cap
H^{s'}(\mathbb{R}) \cap \dot{H}^{-1}(\mathbb{R})$ for $s,s' \geq
2$ and $s \neq s'$, then
\begin{eqnarray*}
u(t) \in C([0,T),H^{s}(\mathbb{R}) \cap \dot{H}^{-1}(\mathbb{R}))\cap
C^{1}([0,T),H^{s-1}(\mathbb{R}))
\end{eqnarray*}
and
\begin{eqnarray*}
u(t) \in C([0,T'),H^{s'}(\mathbb{R}) \cap
\dot{H}^{-1}(\mathbb{R}))\cap C^{1}([0,T'),H^{s'-1}(\mathbb{R}))
\end{eqnarray*}
with the same $T' = T$. See Yin \cite{Yin} for standard arguments.
\end{remark}

By using the local well-posedness result in Lemma \ref{lemma-conservation} and energy
estimates, we obtain a precise blow-up scenario of the solutions to the Cauchy problem (\ref{Cauchy}).

\begin{lemma}
Let $u_{0} \in H^2(\mathbb{R}) \cap \dot{H}^{-1}(\mathbb{R})$ and $u(t)$ be a solution
of the Cauchy problem (\ref{Cauchy}) in Lemma \ref{lemma-conservation}. The solution blows
up in a finite time $T \in (0,\infty)$ in the sense of $\lim_{t \uparrow T} \| u(\cdot,t) \|_{H^2} = \infty$
if and only if
$$
\lim_{t \uparrow T} \sup_{x \in \mathbb{R}} u(x,t) u_x(x,t) =  + \infty.
$$
\label{lemma-blow-up}
\end{lemma}

\begin{proof}
We only need to prove the necessary condition, since the singularity in 
$u(x,t) u_x(x,t)$ as $t \uparrow T$ implies the singularity in $\| u(\cdot,t) \|_{H^2}$ 
as $t \uparrow T$. Assume a finite maximal existence time $T \in (0,\infty)$ and suppose, 
by the contradiction, that there is $M > 0$ such that
\begin{equation}
\label{bound-der} \sup_{x \in \mathbb{R}} u(x,t) u_x(x,t) \le M <
\infty, \quad \forall t \in [0,T).
\end{equation}
Applying density arguments, we approximate the initial value $u_0
\in H^2(\mathbb{R})$ by functions $u_0^n \in H^3(\mathbb{R})$, $n
\geq 1$, so that $\lim_{n \to \infty} u_0^n = u_0$. Furthermore,
write $u^n(t)$ for the solution of the Cauchy problem (\ref{Cauchy})
with initial data $u_0^n$. Using the regularity result proved in
Lemma \ref{lemma-conservation} for $s = 3$, it follows from
Sobolev's embedding that, if $u^n(t) \in C([0,T),H^3(\mathbb{R})\cap
\dot{H}^{-1}(\mathbb{R}))$, then $u^n(x,t)$ is a twice
continuously differentiable function of $x$ on $\mathbb{R}$ for
any $t \in [0,T)$.  It is then deduced from the short-pulse
equation (\ref{short-pulse}) that
$$
\frac{d}{dt} \int_{\mathbb{R}}(u^n_x)^2 dx = \int_{\mathbb{R}} u^n
(u^n_x)^3 dx \le M \int_{\mathbb{R}} (u^n_x)^2 dx
$$
and
$$
\frac{d}{dt} \int_{\mathbb{R}} (u^n_{xx})^2 dx = 5
\int_{\mathbb{R}} u^n u^n_x (u^n_{xx})^2 dx \le 5 M
\int_{\mathbb{R}} (u^n_{xx})^2 dx.
$$
The Gronwall inequality yields for all $t \in [0,T)$,
$$
\|u^n_x(\cdot,t) \|_{L^2} \le \| (u_0^n)' \|_{L^2} e^{\frac{M}{2} t}, \quad 
\|u^n_{xx}(\cdot,t) \|_{L^2} \le \| (u_0^n)'' \|_{L^2} e^{\frac{5}{2} Mt}.
$$
Since $ \|u^n_0\|_{H^2} $ converges to $\|u_0\|_{H^2}$ as $n \to
\infty$, we infer from  the continuous dependence of the local
solution $u(t)$ on initial data $u_0$ that $\| u(\cdot,t) \|_{H^2}$ remains
bounded on $[0,T)$ for the solution $u(t)$ in Lemma
\ref{lemma-conservation}. Therefore, the contradiction is obtained and 
either $T$ is not a maximal existence time or the bound (\ref{bound-der}) is not valid on
$[0,T)$.
\end{proof}

We also show that the blow-up of Lemma \ref{lemma-blow-up} is the wave breaking in
the sense of condition (\ref{def-wave-breaking}). In other words,
both $\| u(\cdot,t) \|_{L^{\infty}}$ and $\| \partial_x^{-1} u(\cdot,t) \|_{L^{\infty}}$ 
are uniformly bounded for all $t \in [0,T)$.

\begin{lemma}
Let $u_{0} \in H^2(\mathbb{R}) \cap \dot{H}^{-1}(\mathbb{R})$ and
$T > 0$ be the maximal existence time of the solution $u(x,t)$ in
Lemma \ref{lemma-conservation}. Then,
\begin{equation}
\label{bounds-u-g} \| u(\cdot,t) \|_{L^{\infty}}
\leq F_1, \quad \| \partial_x^{-1}  u(\cdot,t) \|_{L^{\infty}} \leq F_0, \quad
t \in [0,T),
\end{equation}
where
\begin{eqnarray}
\label{bound-F1} F_1 & := & \frac{1}{\sqrt{2}} \left( E_1^2 + \left( 8 E_0 E_1 + E_1^4 \right)^{1/2} \right)^{1/2}, \\ \label{bound-F0} F_0 & := & \frac{1}{\sqrt{2}} \left( E_0 + E_{-1} +
\frac{1}{12} E_0 F_1^2 \right)^{1/2}.
\end{eqnarray}
\label{lemma-bounds-F0}
\end{lemma}

\begin{proof}
For all $t \in [0,T)$ and the solution $u(x,t)$, we have
\begin{eqnarray*}
u^2(x,t) & = & \left| \int_{-\infty}^{x} u u_x dx - \int_x^{\infty} u u_x dx \right| \\
& \leq & \int_{\mathbb{R}} \frac{|u| |u_x|}{\sqrt{1 + \sqrt{1 + u_x^2}}} \sqrt{1 + \sqrt{1 + u_x^2}} dx \\
& \leq & E_1^{1/2} \left( \int_{\mathbb{R}} u^2 (2 + \sqrt{1 + u_x^2} - 1) dx \right)^{1/2} \\
& \leq & E_1^{1/2} \left( 2 E_0 + E_1 \| u(\cdot,t) \|^2_{L^{\infty}} \right)^{1/2}.
\end{eqnarray*}
As a result, we obtain
$$
\| u(\cdot,t) \|^4_{L^{\infty}} \leq 2 E_0 E_1 + E_1^2 \| u(\cdot,t) \|^2_{L^{\infty}},
$$
so that bound (\ref{bound-F1}) is found from the quadratic equation on 
$\| u(\cdot,t) \|^2_{L^{\infty}}$. On the other hand, we have
\begin{eqnarray*}
\| \partial_x^{-1} u(\cdot,t) \|^2_{H^1} & = & \| u(\cdot,t) \|^2_{L^2} + \|
\partial_x^{-1} u(\cdot,t) \|^2_{L^2} \\
& = & E_0 + E_{-1} + \frac{1}{12} \|
u(\cdot,t) \|^4_{L^4} \\ & \leq & E_0 + E_{-1} + \frac{1}{12} E_0 \|
u(\cdot,t) \|^2_{L^{\infty}}.
\end{eqnarray*}
Using the Sobolev inequality $\| \partial_x^{-1} u \|_{L^{\infty}} \leq \frac{1}{\sqrt{2}} \|
\partial_x^{-1} u \|_{H^1}$ and bound (\ref{bound-F1}), we obtain bound (\ref{bound-F0}).
\end{proof}

Let us introduce a continuous family of characteristics for
solutions of the Cauchy problem (\ref{Cauchy}). Let $\xi \in
\mathbb{R}$, $t \in [0,T)$, where $T$ is the maximal existence
time in Lemma \ref{lemma-conservation}, and denote
\begin{equation}
\label{characteristics-1}
x = X(\xi,t), \quad u(x,t) = U(\xi,t), \quad \partial_x^{-1} u(x,t) = G(\xi,t),
\end{equation}
so that
\begin{equation}
\label{characteristics-2} \left\{ \begin{array}{l} \dot{X}(t) =
-\frac{1}{2} U^2,
\\ X(0) = \xi, \end{array} \right. \quad \left\{ \begin{array}{l}
\dot{U}(t) = G, \\ U(0) = u_0(\xi), \end{array} \right.
\end{equation}
where dots denote derivatives with respect to time $t$ on a
particular characteristics $x = X(\xi,t)$ for a fixed $\xi \in
\mathbb{R}$.  Applying classical results in the theory of ordinary
differential equations, we obtain two useful results
on the solutions of the initial-value problem
(\ref{characteristics-2}). Conserved quantities
 $E_{-1}$ and $E_0$ of the Cauchy
problem (\ref{Cauchy}) are used to control values of $U$ and $G$ on the
family of characteristics.

\begin{lemma}
Let $u_{0}\in H^2(\mathbb{R}) \cap \dot{H}^{-1}(\mathbb{R})$ and
$T>0$ be the maximal existence time of the solution $u(t)$ in
Lemma \ref{lemma-conservation}. Then there exists a unique
solution $X(\xi,t) \in C^{1}(\mathbb{R}\times [0,T))$ to the
initial-value problem (\ref{characteristics-2}). Moreover, the map
$X(\cdot,t) : \mathbb{R} \mapsto \mathbb{R}$ is an increasing
diffeomorphism for any $t \in [0,T)$ with
$$
\partial_{\xi} X(\xi,t) = \exp \left ( -\int_{0}^{t} u(X(\xi,s),s)
u_x(X(\xi,s),s) ds \right ) > 0, \;\; t \in [0,T), \;\;
\xi \in \mathbb{R}.
$$
\label{lemma-characteristics}
\end{lemma}

\begin{proof}
Existence and uniqueness of $X(\xi,t) \in C^{1}(\mathbb{R},[0,T))$ follows from the integral equation
$$
X(\xi,t) = \xi - \frac{1}{2} \int_0^t U^2(\xi,s) ds, \quad t \in
[0,T), \quad \xi \in \mathbb{R},
$$
since $U(\xi,t) \in C(\mathbb{R},[0,T))$ for the solution
$u(t)$ in Lemma \ref{lemma-conservation}. Using the chain rule, we
obtain
$$
\partial_{\xi} \dot{X}(\xi,t) = -W(\xi,t) \partial_{\xi} X(\xi,t)  \quad \Rightarrow \quad
\partial_{\xi} X(\xi,t) = \exp\left(-\int_0^t W(\xi,s) ds\right),
$$
where $W(\xi,t) = u(X(\xi,t),t) u_x(X(\xi,t),t) \in
C(\mathbb{R},[0,T))$. Therefore, $\partial_{\xi} X(\xi,t) > 0$
for all $t \in [0,T)$ and $\xi \in \mathbb{R}$.
\end{proof}

Let
$$
V(\xi,t) = u_x(X(\xi,t),t), \quad W(\xi,t) = u(X(\xi,t),t) u_x(X(\xi,t),t) \equiv U(\xi,t) V(\xi,t)
$$
and compute their rate of changes along the family of characteristics
\begin{eqnarray}
\label{diff-system} \left\{ \begin{array}{lcl}
\dot{V} & = & V W + U, \\
\dot{W} & = & W^2 + V G + U^2.
\end{array} \right.
\end{eqnarray}
Let $F_0, F_1 > 0$ be fixed in terms of conserved quantities $E_{-1}$, $E_0$, and $E_1$ as in Lemma \ref{lemma-bounds-F0} and assume that there exists $\xi_0 \in \mathbb{R}$ such that
$W(\xi_0,0) > 0$ and
\begin{eqnarray*}
& \mbox{\rm either} & \quad |V(\xi_0,0)| > \left( \frac{F_1^2}{4 F_0} \right)^{1/3}, \quad
|V(\xi_0,0)| W(\xi_0,0) > F_1 + \left( 2 F_0 |V(\xi_0,0)|^3 - \frac{1}{2} F_1^2 \right)^{1/2}, \\
& \mbox{\rm or} & \quad |V(\xi_0,0)| \leq \left( \frac{F_1^2}{4 F_0} \right)^{1/3}, \quad
|V(\xi_0,0)| W(\xi_0,0) > F_1.
\end{eqnarray*}
Because of the invariance of system (\ref{diff-system}) with respect to
$$
G \to -G, \quad U \to -U, \quad V \to -V, \quad W \to W,
$$
it is sufficient to consider the case with $V(\xi_0,0) > 0$.
We will prove that, under the conditions above, $V(\xi_0,t)$ and $W(\xi_0,t)$ remain
positive and monotonically increasing functions for all $t > 0$,
for which they are bounded, so that $V(\xi_0,t)$ and $W(\xi_0,t)$ satisfy
the apriori differential estimates
\begin{eqnarray}
\label{diff-system-apriori}
\left\{ \begin{array}{lcl}
\dot{V} & \geq & V W - F_1, \\
\dot{W} & \geq & W^2 - V F_0.
\end{array} \right.
\end{eqnarray}
In what follows, we use $V(t)$ and $W(t)$ instead of $V(\xi_0,t)$ and
$W(\xi_0,t)$ for a particular $\xi_0 \in \mathbb{R}$. The following lemma
establishes sufficient conditions on the initial point
$(V(0),W(0))$ that ensure that a lower solution satisfying
\begin{eqnarray}
\label{diff-system-minor}
\left\{ \begin{array}{lcl}
\dot{V} & = & V W - F_1, \\
\dot{W} & = & W^2 - V F_0,
\end{array} \right.
\end{eqnarray}
goes to infinity in a finite time.

\begin{lemma}
\label{lemma-trajectories}
Assume that the initial data for system (\ref{diff-system-minor})
satisfy
\begin{eqnarray}
\label{suff-cond}
& \mbox{\rm either} & \quad V(0) > \left( \frac{F_1^2}{4 F_0} \right)^{1/3}, \quad
V(0) W(0) > F_1 + \left( 2 F_0 V^3(0) - \frac{1}{2} F_1^2 \right)^{1/2}, \\
\label{suff-cond-a}
& \mbox{\rm or} & \quad 0 < V(0) \leq \left( \frac{F_1^2}{4 F_0} \right)^{1/3}, \quad
V(0) W(0) > F_1.
\end{eqnarray}
Then the trajectory of system (\ref{diff-system-minor}) blows up in a finite time
$t_* \in (0,\infty)$ such that $V(t)$ and $W(t)$ are positive and monotonically increasing
for all $t \in [0,t_*)$ and there is $C > 0$ such that
\begin{equation}
\label{limit-cond}
\lim_{t \uparrow t_*} V(t) = \infty, \quad
\lim_{t \uparrow t_*} W(t) = \infty, \quad \mbox{\rm and} \quad \lim_{t \uparrow t_*} (t_* - t) V(t) = C.
\end{equation}
Moreover, $t_*$ is bounded by
\begin{equation}
\label{bound-on-time}
t_* \leq \frac{V(0)}{\min\left\{ \dot{V}(0), \left(\dot{V}^2(0) - 2 F_0 V^3(0) + \frac{1}{2} F_1^2\right)^{1/2} \right\}}.
\end{equation}
\end{lemma}

\begin{proof}
Let us first consider the homogeneous version of system (\ref{diff-system-minor}) for $F_1 = 0$,
that is
\begin{equation}
\label{simple-system}
\left\{ \begin{array}{lcl}
\dot{V} & = & V W, \\
\dot{W} & = & W^2 - V F_0.
\end{array} \right.
\end{equation}
Of course, $F_1$ is never zero, otherwise $E_1 = 0$. This case is used merely for illustration,
since explicit solutions can be obtained for $F_1 = 0$, whereas qualitative analysis has to be
developed for $F_1 \neq 0$. System (\ref{simple-system}) is integrable since
$$
W = \frac{\dot{V}}{V} \quad \Rightarrow \quad \frac{d}{dt} \left( \frac{\dot{V}}{V^2} \right) = -F_0
\quad \Rightarrow \quad \dot{V} = V^2 (C - F_0 t),
$$
where $C = W(0)/V(0)$.
Integrating the last equation for $V(t)$, we obtain
the explicit solution of the truncated system,
$$
V(t) = \frac{V(0)}{1 - C V(0) t + \frac{1}{2} F_0 V(0) t^2}, \quad W(t) = (C - F_0 t) V(t).
$$
The solution reaches infinity in a finite time $t_* \in (0,\infty)$
if $V(0) > 0$ and $C^2 V(0) > 2 F_0$. Note that these conditions coincide with
condition (\ref{suff-cond}) for $F_1 = 0$. Also note that $t_*$ is the first positive root of
$1 - C V(0) t + \frac{1}{2} F_0 V(0) t^2 = 0$ so that $U(t_*) = C - F_0 t_* > 0$ and 
$$
t_* = \frac{W(0) - \sqrt{W^2(0) - 2 F_0 V(0)}}{F_0 V(0)} = \frac{2 V(0)}{\dot{V}(0) + \sqrt{\dot{V}^2(0) - 2 F_0 V^3(0)}} \leq \frac{V(0)}{\sqrt{\dot{V}^2(0) - 2 F_0 V^3(0)}}.
$$

Consider now the full system (\ref{diff-system-minor}). Let
$$
V = \frac{1}{x}, \quad W = \frac{y}{x}
$$
and rewrite the system in the form
\begin{equation}
\label{diff-system-x}
\left\{ \begin{array}{lcl}
\dot{x} & = & F_1 x^2 - y, \\
\dot{y} & = & F_1 x y - F_0.
\end{array} \right.
\end{equation}
Expressing $y$ from the system,
we can rewrite it in the scalar form
\begin{equation}\label{scalar-ode}
y = F_1 x^2 - \dot{x} \quad \Rightarrow \quad \ddot{x} = f(x) + 3 F_1 x \dot{x},
\end{equation}
where $f(x) = F_0 - F_1^2 x^3$. The only critical point of equation (\ref{scalar-ode}) is $(x,\dot{x}) = (x_0,0)$ where $x_0 = (F_0/F_1^2)^{1/3}$ is the root of $f(x)$.

We need to show that there is a domain $D \subset \R^2$ in the phase plane $(x,\dot{x})$,
so that all initial data in $D$ generate trajectories in $D$ that cross the vertical line $x = 0$
in a finite time. To do so, we construct a Lyapunov function for
equation (\ref{scalar-ode}) in the form
$$
E(x,\dot{x}) = \frac{1}{2} \dot{x}^2 - F_0 x + \frac{1}{4} F_1^2 x^4.
$$
The function $E(x,\dot{x})$ has a global minimum at $(x_0,0)$. For any solution $x(t)$, we have
$$
\frac{d}{dt} E(x,\dot{x}) = 3 F_1 x \dot{x}^2 > 0 \quad \mbox{\rm for} \quad x > 0.
$$
The zero level of the Lyapunov function $E(x,\dot{x}) = 0$
passes through the points $(0,0)$ and $(x_*,0)$, where $x_* = 4^{1/3} x_0  > x_0$
(see Figure \ref{fig:PhasePlane}). It is clear that $E(x,\dot{x}) > 0$ in the domain
$$
D = \left\{ (x,\dot{x}): \quad x > 0, \quad \dot{x} < \sigma(x) \right\},
$$
where
$$
\sigma(x) = \left\{ \begin{array}{cl}
-\left( 2 F_0 x - \frac{1}{2} F_1^2 x^4 \right)^{1/2}, \quad & x \in (0,x_*), \\
0, \quad & x \in [x_*,\infty). \end{array} \right.
$$
We note that the condition $(x,\dot{x}) \in D$ is equivalent to the condition
$$
\left\{ \begin{array}{ll} 0 < x < x_*, \quad & y > F_1 x^2 + \left( 2 F_0 x - \frac{1}{2} F_1^2 x^4 \right)^{1/2}, \\
x > x_*, \quad & y > F_1 x^2, \end{array} \right.
$$
which is nothing but the set of conditions (\ref{suff-cond}) and
(\ref{suff-cond-a}) at $t = 0$. By continuity, if $(x,\dot{x}) \in D$ at $t = 0$,
then $(x,\dot{x})$ remains in $D$ for some time $t > 0$.

No critical points of system (\ref{scalar-ode}) are located in $D$ and $x(t)$
is decreasing function for any $t > 0$ as long as the trajectory stays in $D$.
Recall that $E(x,\dot{x}) > 0$ and $\frac{d}{dt} E(x,\dot{x}) > 0$ for any $(x,\dot{x}) \in D$.
A trajectory in $D$ can not cross $\dot{x} = \sigma(x)$ because
$E(x,\sigma(x)) = 0$ for $x \in (0,x_*)$ and
$\ddot{x} = f(x) + 3 F_1 x \dot{x} < 0$ for $x > x_*$ and $\dot{x} < 0$.
Therefore, the trajectory either reaches $x = 0$ in a finite time $t_* \in (0,\infty)$
or escapes to $\dot{x} = -\infty$ for $x > 0$. To eliminate the last possibility, we note
that
$$
\frac{d}{dt} \left( \dot{x} - \frac{3 F_1}{2} x^2 \right) = f(x),
$$
so that
$$
\dot{x}(t) \geq \dot{x}(0) - \frac{3F_1}{2} x^2(0) + t f(x(0)) > -\infty
$$
for any finite time interval. Moreover, $\dot{x}$ is bounded from zero in $D$ by
the level curve $E(x,\dot{x}) = E(x(0),\dot{x}(0))$, which is a convex curve in $D$. Therefore,
$$
\dot{x}(t) \leq \max\{ \dot{x}(0),\rho \}, \quad t > 0, \quad \mbox{\rm as long as} \quad x > 0,
$$
where $\rho < 0$ is uniquely found from $E(0,\rho) = E(x(0),\dot{x}(0))$, 
that is from the point of intersection of the level curve of $E(x,\dot{x}) = E(x(0),\dot{x}(0))$ 
with the negative $\dot{x}$-axis. Therefore,
$$
x(t) \leq x(0) + t \max\{ \dot{x}(0),\rho \},
$$
so that $x(t)$ reaches $x = 0$ in a finite time $t_* \in (0,\infty)$
for any trajectory in $D$. Moreover, finding $\rho$ explicitly gives the bound on the
blow-up time
$$
\left\{ \begin{array}{ll} 0 < x < x_*, \quad & t_* \leq \frac{x(0)}{(\dot{x}^2(0) - 2 F_0 x(0) + \frac{1}{2} F_1^2 x^4(0))^{1/2}}, \\ \vspace{0cm} \\
x > x_*, \quad & t_* \leq \frac{x(0)}{|\dot{x}(0)|}, \end{array} \right.
$$
which becomes bound (\ref{bound-on-time}) after the return back to variable $V(t)$.

\begin{figure}
\begin{center}
\includegraphics[width=0.45\textwidth]{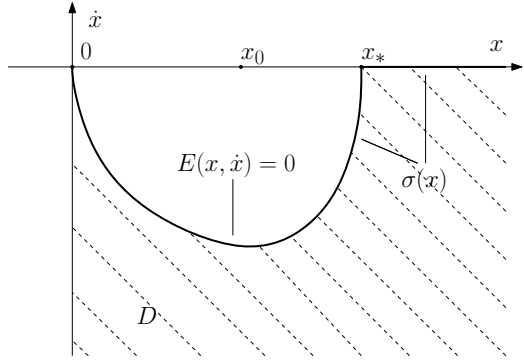}
\end{center}
\caption{Domain $D$ in the phase plane $(x,\dot{x})$ of equation (\ref{scalar-ode}). }
\label{fig:PhasePlane}
\end{figure}

Since $V = x^{-1}$, $W = y x^{-1}$, and $y = F_1 x^2 - \dot{x} > 0$, we have
$\lim_{t \uparrow t_*} V(t) = \infty$ and $\lim_{t \uparrow t_*} W(t) = \infty$
for any trajectories in $D$. Since $\dot{x} < 0$ for the trajectory in $D$, we also have
$x(t) \sim (t_* - t)$ as $t \to t_*$ so that there exists $C > 0$ such that
$\lim_{t \uparrow t_*} (t_* - t) V(t) = C$. It remains to show that
$V(t)$ and $W(t)$ are monotonically increasing functions on $[0,t_*)$. To do so,
we write
\begin{eqnarray*}
\dot{V} & = & V W - F_1 = -\frac{\dot{x}}{x^2}, \\
\dot{W} & = & W^2 - F_0 V = \frac{g(x,\dot{x})}{x^2},
\end{eqnarray*}
where
$$
g(x,\dot{x}) = \dot{x}^2 - 2 F_1 x^2 \dot{x} - x f(x).
$$
For any trajectory in $D$, $\dot{x}(t) < 0$ so that $\dot{V}(t) > 0$. Furthermore,
since $g(x,\dot{x})$ is zero at a curve outside the domain $D$, because
$$
g(x,\sigma(x)) = F_0 x + \frac{1}{2} F_1^2 x^4 + 2 F_1 x^2 |\sigma(x)| > 0, \quad x \in (0,x_*],
$$
then $\dot{W}(t) > 0$ for any trajectory in $D$.
\end{proof}

Recalling that $W = U V$ in system (\ref{diff-system}) and using
the bound (\ref{bounds-u-g}), we obtain the upper bound for any
solution at the family of characteristics
\begin{equation}
\label{super-solution} \dot{V} = U V^2 + U \leq F_1 V^2 + F_1.
\end{equation}
We can now show that any upper solution with $V(0) > 0$ goes to
infinity in a finite time.

\begin{lemma}
\label{lemma-supersolution} Consider
\begin{equation}
\label{separable-equation} \dot{V}(t) = F_1 V^2 + F_1
\end{equation}
with $V(0) > 0$. There exists $t_* \in (0,\infty)$ such that
$V(t)$ is positive, monotonically increasing for all $t \in
[0,t_*)$ and there exists $C > 0$ such that
$$
\lim_{t \uparrow t_*} (t_* - t) V(t) = C.
$$
Moreover, $t_* \leq 1/W(0)$.
\end{lemma}

\begin{proof}
Since $\dot{V} > 0$ for any $V \in \R$, $V(t)$ is monotonically increasing function.
To show that $V(t)$ reaches $\infty$ in a finite time, one can integrate the
separable equation (\ref{separable-equation}) explicitly and obtain
$$
V(t) = \tan( \arctan V(0) + F_1 t),
$$
so that
$$
t_* = \frac{\pi/2 - \arctan V(0)}{F_1} \leq \frac{1}{F_1 V(0)} \leq \frac{1}{W(0)},
$$
since $\sup_{x \in \R_+} x \cot^{-1}(x) \leq 1$.
\end{proof}

Applying results of Lemmas \ref{lemma-trajectories} and
\ref{lemma-supersolution}, we conclude the proof of Theorem
\ref{theorem-wave-breaking}.

\begin{proof1}{\em of Theorem \ref{theorem-wave-breaking}.}
Let $(V,W)$ satisfy system (\ref{diff-system}) corresponding to
the characteristics with $\xi_0$. Let
$(\underline{V},\underline{W})$ be the lower solution of system
(\ref{diff-system-minor}) in Lemma \ref{lemma-trajectories} with
$\underline{V}(0) = V(0)$ and $\underline{W}(0) = W(0)$. Let $\overline{V}$ be the
upper solution of equation (\ref{separable-equation}) in Lemma
\ref{lemma-supersolution} with $\overline{V}(0) = V(0)$. Let $\underline{t}_*$
be the blow-up time of the lower solution and $\overline{t}_*$ be the
blow-up time of the upper solution.

The upper bound for the solution of system (\ref{diff-system}) follows from
the comparison principle for the differential equations since
$$
|V W + U| = (V^2 + 1) |U| \leq (V^2 + 1) F_1
$$
which implies that $V(t) \leq \overline{V}(t)$ for all $t \in
[0,\overline{t}_*)$, for which $\overline{V}(t)$ remains bounded.

To obtain the lower bound, we note that
$$
V \geq \underline{V}, \; W \geq \underline{W} \quad \Rightarrow
\quad \left\{
\begin{array}{l}
\dot{V} \geq VW - F_1 \geq \underline{V} \underline{W} - F_1 = \dot{\underline{V}}, \\
\dot{W} + F_0 V \geq W^2 \geq \underline{W}^2 =
\dot{\underline{W}} + F_0 \underline{V}.
\end{array} \right.
$$
Let ${\bf V} = [V,W]^T$, $\underline{\bf V} =
[\underline{V},\underline{W}]^T$, and
$$
L = \left[ \begin{array}{cc} 0 & 0 \\ F_0 & 0 \end{array} \right]
$$
be a nilpotent matrix of order one, so that $e^{t L} = I + t L$. Thus, we write
$$
V \geq \underline{V}, \; W \geq \underline{W} \quad \Rightarrow
\quad \frac{d}{dt} \left( e^{t L} {\bf V} \right) \geq
\frac{d}{dt} \left( e^{t L} \underline{\bf V} \right).
$$
Integrating this equation in $t > 0$, we infer that
$$
e^{tL} {\bf V}(t) \geq e^{tL} \underline{\bf V}(t),
$$
and since $e^{tL}$ is invertible for any $t \in \R$, we conclude
that
$$
V(t) \geq \underline{V}(t), \quad W(t) \geq \underline{W}(t),
$$
for all $t \in [0,T) \subset [0,\underline{t}_*)$, for which $(V,W)$ remain finite.
Therefore,
$(V,W)$ become infinite as $t \uparrow T$ and $T \in [\overline{t}_*,\underline{t}_*]$.
\end{proof1}

\begin{remark}
The bounds on $\underline{t}_*$ and $\overline{t}_*$ in Lemmas \ref{lemma-trajectories} and
\ref{lemma-supersolution} are inconclusive to compare $T$ with the exact time of blow-up 
$T_0 := \frac{1}{W(0)}$ along a particular characteristic of the 
dispersionless advection equation (\ref{simple-wave}).
\end{remark}

\section{Wave breaking in a periodic domain}

Consider now the Cauchy problem for the short-pulse equation
(\ref{short-pulse}) in a periodic domain
\begin{equation}
\label{cauchy-periodic}
\left\{ \begin{array}{ll} u_t = \frac{1}{2} u^2 u_x + \partial_x^{-1} u, \quad & x \in \mathbb{S}, \;\; t > 0,\\
u(x,0) = u_0(x), \quad & x \in \mathbb{S}, \end{array} \right.
\end{equation}
where $\mathbb{S}$ is a unit circle equipped with periodic boundary conditions and
$\partial_x^{-1}$ is the mean-zero anti-derivative in the form
$$
\partial_x^{-1} u := \int_0^x u(x',t) dx' - \int_{\mathbb{S}} \int_0^x u(x',t) dx' dx.
$$

Local well-posedness and useful conserved quantities for the
Cauchy problem (\ref{cauchy-periodic}) in a periodic domain are
obtained in the following lemma.

\begin{lemma}
Assume that $u_{0} \in H^s(\mathbb{S})$, $s \geq 2$ and
$\int_{\mathbb{S}} u_0(x) \, dx = 0$. Then  there exist a maximal
time $T > 0$ such that the Cauchy problem (\ref{cauchy-periodic})
admits a unique solution
$$
u(t) \in C([0,T),H^s(\mathbb{S})) \cap C^{1}([0,T),H^{s-1}(\mathbb{S}))
$$
satisfying $u(x,0) = u_0(x)$ and $\int_{\mathbb{S}} u(x,t) dx  = 0$ for
all $t \in [0,T)$. Moreover, the solution $u(t)$ depends
continuously on the initial data $u_0$ and the quantities
$$
E_0 = \int_{\mathbb{S}} u^2 dx, \quad
E_1 = \int_{\mathbb{S}} \sqrt {1 + u^2_x} dx
$$
are constant for all $t \in [0,T)$. \label{lemma-Wayne}
\end{lemma}

\begin{proof}
Existence of the solution $u(x,t)$ and continuous dependence on
$u_0$ is proved on $\mathbb{S}$ similarly to what is done
in Lemma \ref{lemma-conservation} on $\R$. To prove the zero-mass
constraint, we note
$$
u_t(t) \in C((0,T),H^{s-1}(\mathbb{S})), \quad u^2 u_x(t) \in C([0,T),H^{s-1}(\mathbb{S})),
$$
so that for all $t \in (0,T)$, we have
$$
\int_{\mathbb{S}} u(x,t) \, dx = \int_{\mathbb{S}} u_{tx} \, dx
+ \frac{1}{2} \int_{\mathbb{S}} (u^2 u_x)_x \, dx
= 0.
$$
Initial values of $E_0$ and $E_1$ are bounded if $u_0 \in
H^s(\mathbb{S})$, $s \geq 2$. Conservation of $E_0$ and $E_1$ on
$[0,T)$ follows from the balance equations
\begin{eqnarray*}
\left( u^2 \right)_t & = & \left( (\partial_x^{-1} u)^2 + \frac{1}{4} u^4 \right)_x, \\
\left( \sqrt{1 + u_x^2} \right)_t & = & \left( \frac{1}{2} u^2 \sqrt{1 + u_x^2} \right)_x,
\end{eqnarray*}
thanks to the continuity and the periodic boundary conditions for
$\partial_x^{-1} u(t) \in C((0,T),H^{s+1}(\mathbb{S}))$, $u(t) \in
C((0,T),H^s(\mathbb{S}))$, and $u_x(t) \in
C((0,T),H^{s-1}(\mathbb{S}))$ in $x$ on $\mathbb{S}$ if $s \geq
2$.
\end{proof}

\begin{remark} The assumption $\int_{\mathbb{S}} u_0(x) \, dx = 0$ on the initial data $u_0$
in the periodic domain $\mathbb{S}$ is necessary as it follows from the following apriori estimate
$$
\left | \int_{\mathbb{S}} u(x,t) \, dx - \int_{\mathbb{S}} u_0 (x)
\, dx \right | \le \| u(t) - u_0 \|_{L^2(\mathbb{S})}, \quad
\forall t \in (0, T).
$$
Note that $\int_{\mathbb{S}} u(x,t) \, dx = 0$, for all $t \in (0,
T)$ and $u(t) \in C([0,T), H^s(\mathbb{S}))$, $s \geq 2$. Hence
the above estimate implies that $ \int_{\mathbb{S}} u_0(x) \, dx =
0$. Note that no zero-mass constraint is necessary on an infinite
line in Theorem \ref{theorem-wayne}.
\end{remark}

The blow-up scenario for the solutions to the Cauchy problem
(\ref{cauchy-periodic}) coincides with the one in Lemma
\ref{lemma-blow-up} after the change of $\mathbb{R}$ by
$\mathbb{S}$. The main result of this section is the proof of the
finite-time wave breaking in a periodic domain, according to the
following theorem.

\begin{theorem}
\label{theorem-breaking-periodic}
Let $u_{0} \in H^2(\mathbb{S})$ and $\int_{\mathbb{S}}u_0(x) \, dx = 0$.
Assume that there exists $x_0 \in \mathbb{R}$ such that $u_0(x_0) u_0'(x_0) > 0$ and
\begin{eqnarray*}
& \mbox{\rm either} & \quad |u_0'(x_0)| > \left( \frac{E_1^2}{4 E_0^{1/2}} \right)^{1/3}, \quad
|u_0(x_0)| |u_0'(x_0)|^2 > E_1 + \left( 2 E_0^{1/2} |u_0'(x_0)|^3 -
\frac{1}{2} E_1^2 \right)^{1/2}, \\
& \mbox{\rm or} & \quad
|u_0'(x_0)| \leq \left( \frac{E_1^2}{4 E_0^{1/2}} \right)^{1/3}, \quad
|u_0(x_0)| |u_0'(x_0)|^2 > E_1.
\end{eqnarray*}
Then there exists a finite time $T \in (0,\infty)$ such that the
solution $u(t) \in C([0,T),H^2(\mathbb{S}))$ of the Cauchy problem
(\ref{cauchy-periodic}) blows up with the property
$$
\lim_{t \uparrow T} \sup_{x \in \mathbb{S}} u(x,t) u_{x}(x,t) = +
\infty, \quad \mbox{while} \quad \lim_{t \uparrow T} \| u(\cdot,t)
\|_{L^{\infty}} \leq E_1.
$$
\end{theorem}

\begin{proof}
Let $T > 0$ be the maximal time of existence of the solution $u(t)
\in C([0,T),H^2(\mathbb{S}))$ to the Cauchy problem
(\ref{cauchy-periodic}) constructed in Lemma \ref{lemma-Wayne}.
Since $\int_{\mathbb{S}} u(x,t) dx = 0$, for each $t\in [0,T)$
there is a $\xi_{t}\in [0,1]$ such that $u(\xi_{t},t)=0$. Then for
$x\in \mathbb{S}$ and $t \in [0,T)$, we have
\begin{equation*}
|u(x,t)| = \left| \int_{\xi_t}^{x} u_x(x,t) \,dx \right| \leq \int_{\mathbb{S}} |u_x(x,t)| dx \leq E_1.
\end{equation*}
Since $\partial_x^{-1} u(t) \in C([0,T),H^3(\mathbb{S}))$ is the mean-zero periodic function of $x$
for each $t \in [0,T)$, there exists another $\tilde{\xi}_t \in [0,1]$ such that
$\partial_x^{-1} u(\tilde{\xi}_t,t) = 0$. Then for $x\in \mathbb{S}$ and $t \in [0,T)$, we have
\begin{equation*}
|\partial_x^{-1} u(x,t)| = \left| \int_{\tilde{\xi}_t}^{x} u(x,t) \,dx \right| \leq
\int_{\mathbb{S}} |u(x,t)| dx \leq \sqrt{E_0}.
\end{equation*}
Therefore, bounds (\ref{bounds-u-g}) are rewritten with 
$$
F_0 := \sqrt{E_0}, \quad F_1 := E_1
$$
The rest of the proof follows the proof of Theorem \ref{theorem-wave-breaking}.
\end{proof}

\section{Numerical evidence of wave breaking}

The goal of this section is to complement the analytic results by
several examples and numerical computations. More specifically, we
first show that the sufficient condition for wave breaking in
Theorem \ref{theorem-wave-breaking} is not satisfied for the exact
modulated pulse solution to the short-pulse equation which is
known to remain globally bounded in space and time. Then we
consider the interplay between global well-posedness and wave
breaking of Theorems \ref{theorem-wellposedness} and
\ref{theorem-wave-breaking} for a class of decaying data on an
infinite line. Finally, we perform numerical simulations in a
periodic domain for a simple harmonic initial data and thus give
illustrations to the sufficient condition for wave breaking in
Theorem \ref{theorem-breaking-periodic}.

Theorem \ref{theorem-wave-breaking} gives a sufficient condition
for formation of shocks in the short-pulse equation
(\ref{short-pulse}) on the infinite line. Let us show that this
condition is not satisfied for exact modulated pulse solutions
obtained in \cite{Matsuno,SS06}. The simplest one-pulse solution
is given in the parametric form
$$
u(x,t) = U(y,t), \quad x = X(y,t),
$$
where
\begin{equation}
\left\{
\begin{array}{l} \label{pulse1}
U(y,t) = 4 m n \dfrac{m \sin \psi \sinh \phi + n \cos \psi \cosh \phi}{m^2 \sin^2 \psi + n^2 \cosh^2 \phi}, \\[6pt] \\
X(y,t) = y + 2 m n \dfrac{m \sin 2 \psi - n \sinh 2 \phi}{m^2
\sin^2 \psi + n^2 \cosh^2 \phi},
\end{array}
\right. \quad (y,t) \in \R^2,
\end{equation}
$m \in (0,1)$ is an arbitrary parameter, $n = \sqrt{1 - m^2}$, and
$$
\phi = m (y + t), \quad \psi = n (y - t).
$$
The pulse solution enjoys the periodicity property
\begin{equation*}
\left\{
\begin{array}{l}
U(y,t) = U\left(y - \frac{\pi}{m}, t + \frac{\pi}{m}\right), \\[6pt]
X(y,t) = X\left(y - \frac{\pi}{m}, t + \frac{\pi}{m}\right) +
\frac{\pi}{m}
\end{array}
\right. \quad (y,t) \in \mathbb{R}^2.
\end{equation*}
and an exponential decay in any direction transverse to the
anti-diagonal on the $(y,t)$-plane.

\begin{figure}
\begin{center}
\includegraphics[width=0.45\textwidth]{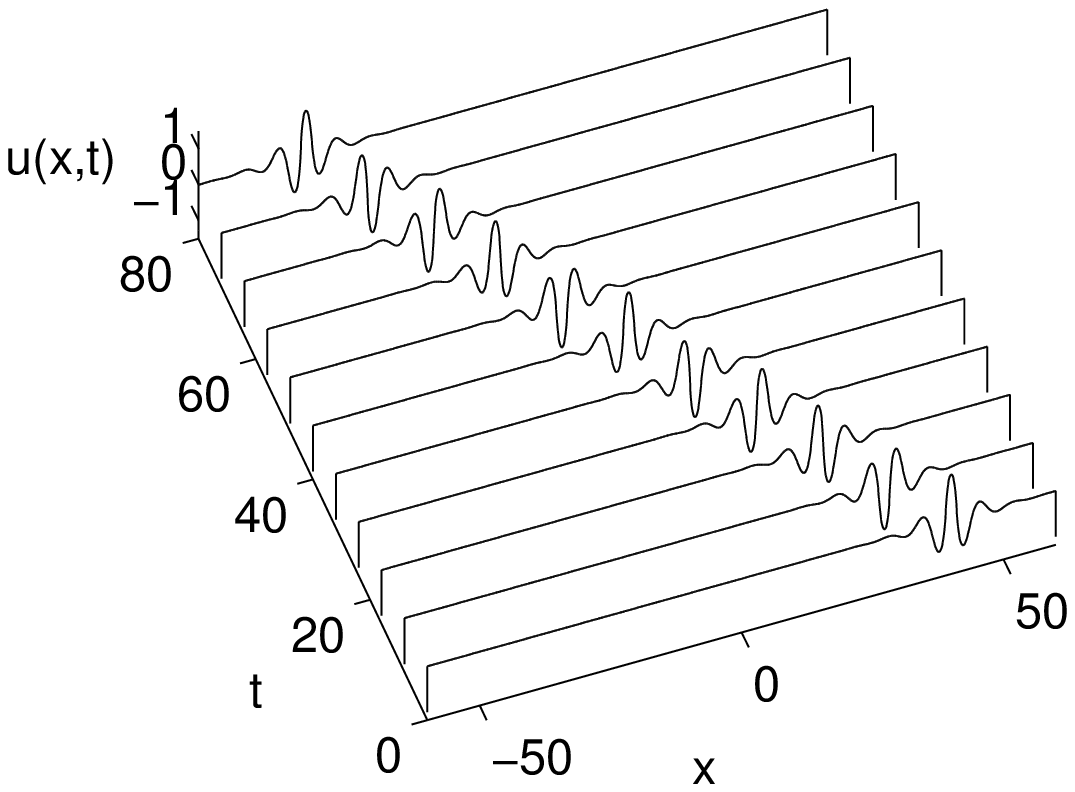}
\includegraphics[width=0.45\textwidth]{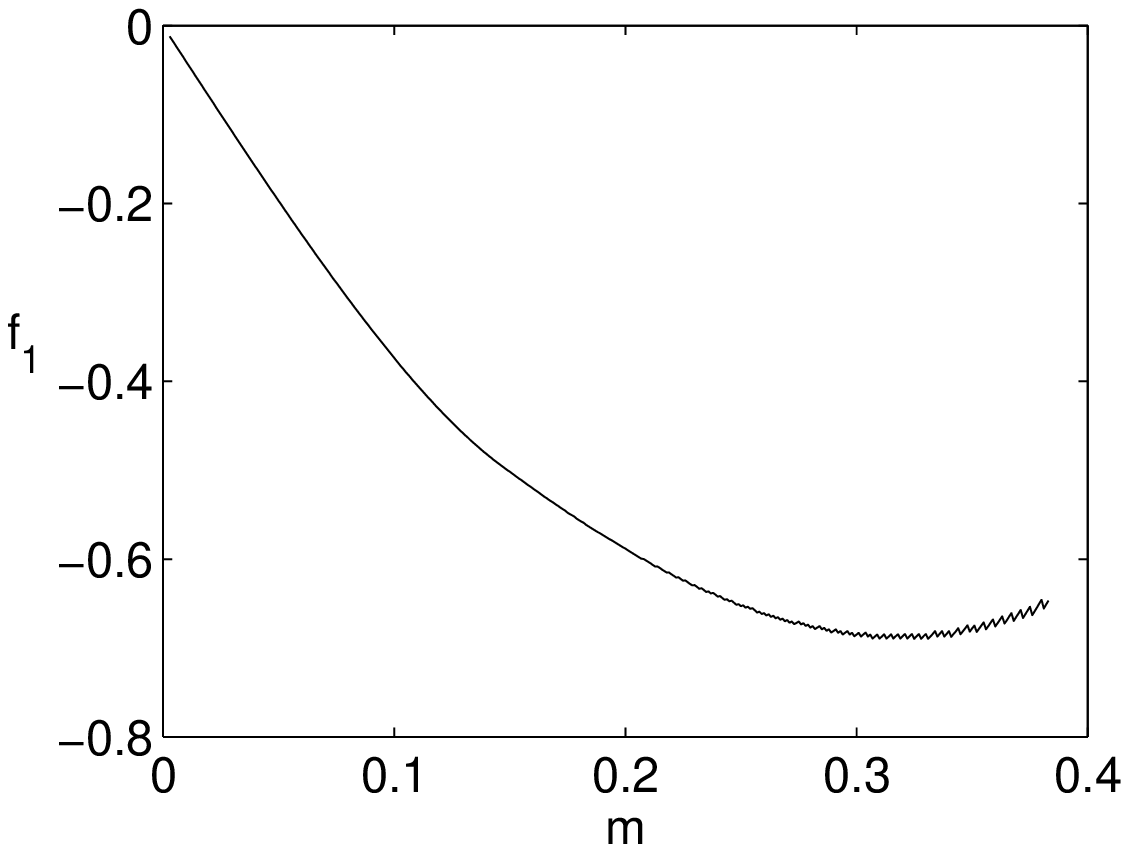}
\end{center}
\caption{The exact modulated pulse solution (\ref{pulse1}) of the
short-pulse equation (\ref{short-pulse}) for $m=0.32$ (left). The
quantity $f_1$ is negative for any $m \in (0,m_{cr})$ (right).}
\label{fig:pulse}
\end{figure}

Since
$$
\frac{\partial X}{\partial y} = 1 - \frac{8 m^2 n^2 \sin^2 \psi \cosh^2 \phi}{(m^2 \sin^2 \psi +
n^2 \cosh^2 \phi)^2} = \cos \left( 4 \arctan \frac{m \sin \psi}{n \cosh \phi} \right),
$$
the function $x = X(y,t)$ is invertible in $y$ for all
$t \in \mathbb{R}$ if
$$
\left| \frac{m \sin \psi}{n \cosh \phi} \right| < \tan\frac{\pi}{8} \quad \Rightarrow
\quad \frac{m}{n} \leq  \tan\frac{\pi}{8},
$$
that is for all $m \in (0,m_{cr})$, where $m_{cr} =
\sin{\frac{\pi}{8}} \approx 0.383$. For these values of $m$, the
pulse solution $u(x,t)$ is analytic in variables
$(x,t) \in \R$, has the space-time periodicity
$$
u(x,t) = u\left(x - \frac{\pi}{m},t + \frac{\pi}{m}\right),  \quad
(x,t) \in \R,
$$
and the exponential decay in the transverse direction to the
anti-diagonal in the $(x,t)$-plane. The graph of a nonsingular
pulse solution for $m = 0.32$ is shown on Figure \ref{fig:pulse}
(left).

\begin{remark}
Coordinate $y$ in the exact solution (\ref{pulse1}) is different from
coordinate $\xi$ in the method of characteristics because $X(y,0) \neq y$.
Nevertheless, $X(y,t)$ and $U(y,t)$ satisfy the same set of equations
$$
\frac{\partial X}{\partial t} = -\frac{1}{2} U^2(y,t), \quad
\frac{\partial U}{\partial t} = \partial_x^{-1} u |_{x = X(y,t)},
$$
so that $\xi$ and $y$ are uniquely related by the representation $\xi = X(y,0)$. If 
$y$ is found as a function of $\xi$, the initial data of the Cauchy problem 
(\ref{Cauchy}) is found from $u_0(\xi) = U(y,0)$.
\end{remark}

Since $u(x,t)$ is analytic in $x \in \R$ for any fixed $t \in \R$
and decays to zero exponentially fast at infinity, it is clear
that $u(\cdot,t) \in H^2(\mathbb{R})$. Furthermore, since
$\partial_x^{-1} u = u_t - \frac{1}{2}u^2 u_x$,
it also follows that $u(\cdot,t) \in \dot{H}^{-1} (\mathbb{R})$.
We compute numerically bounds $F_0$ and $F_1$ using the exact solution (\ref{pulse1})
and new definitions
$$
F_1 := \sup_{t \in \R} \| u(\cdot,t) \|_{L^{\infty}}, \quad F_0 := \sup_{t \in \R}
\|  \partial_x^{-1} u(\cdot,t) \|_{L^{\infty}}.
$$
It follows from Remark \ref{remark-scaling-invariance} that bounds $F_1$ and $F_0$ defined above
preserve the sufficient condition of Theorem \ref{theorem-wave-breaking} with respect to 
the scaling transformation (\ref{scaling-invariance}).

Let us define
\begin{equation}
\label{WB1}
\begin{split}
f_1 & := \sup_{x \in I_1} \left[ |u_0(x)| |u_0'(x)|^2 - F_1 \right], \\
f_2 & := \sup_{x \in I_2} \left[ |u_0(x)| |u_0'(x)|^2 - F_1 - \left( 2F_0
|u_0'(x)|^3 - \frac{1}{2} F_1^2 \right)^{1/2}  \right],
\end{split}
\end{equation}
where
\begin{eqnarray*}
I_1 & = & \left\{ x \in \R \; : \; |u_0'(x)| \leq \left( \dfrac{F_1^2}{4
F_0} \right)^{1/3}, \;\; u_0(x) u_0'(x) > 0, \right\}, \\
I_2 & = & \left\{ x \in \R \; : \; |u_0'(x)| > \left(\dfrac{F_1^2}{4 F_0}
\right)^{1/3}, \;\; u_0(x) u_0'(x) > 0, \right\}.
\end{eqnarray*}
According to Theorem \ref{theorem-wave-breaking}, wave breaking
occurs if either $f_1$ or $f_2$ is positive. For the exact
modulated pulse solution (\ref{pulse1}) at $t = 0$, the numerical
calculations show that the set $I_2$ is empty and the quantity
$f_1$ is strictly negative for any $m \in (0,m_{cr})$, see Figure
\ref{fig:pulse} (right). Therefore, the sufficient condition for
the wave breaking in Theorem \ref{theorem-wave-breaking} is not
satisfied, which corresponds to our understanding that the exact
modulated pulse solutions (\ref{pulse1}) remain bounded for all $(x,t) \in \R^2$.
We note, however, that the sufficient condition for the global
well-posedness in Theorem \ref{theorem-wellposedness} is satisfied
only for pulses with $m \in \left( 0,\frac{1}{32} \right)$, since $2 \sqrt{2
E_1 E_2} = 32 m$. This computation shows that the sufficient
condition of Theorem \ref{theorem-wellposedness} is not sharp.

\begin{figure}
\begin{center}
\includegraphics[width=0.45\textwidth]{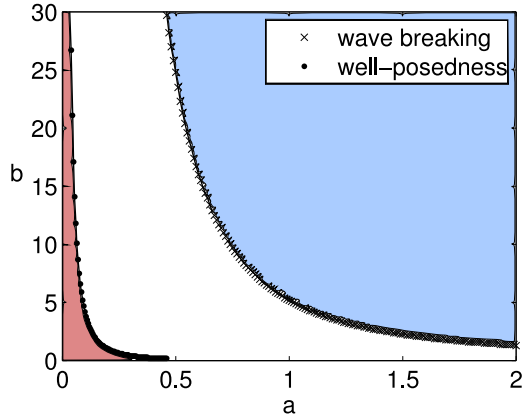}
\end{center}
\caption{Boundaries of the global well-posedness and the wave
breaking in the Cauchy problem (\ref{Cauchy}) with initial data
(\ref{gauss}): the global well-posedness occurs below the lower
curve and the wave breaking occurs above the upper curve.}
\label{fig:gauss}
\end{figure}

Next we compare the sufficient conditions for the global
well-posedness and wave breaking in Theorems
\ref{theorem-wellposedness} and \ref{theorem-wave-breaking} for a
class of initial data
\begin{equation}
u_0(x) = a(1-2 b x^2)e^{-b x^2}, \quad a>0, \quad b > 0,
\label{gauss}
\end{equation}
where parameters $a$ and $b$ determine the amplitude and steepness
of $u_0$. Note that the zero-mass constraint
(\ref{mass-conservation}) is satisfied by $u_0$ and it is clear
that $u_0 \in H^2(\R) \cap \dot{H}^{-1}(\R)$. The
conserved quantities $E_{-1}$ and $E_0$ can be computed
analytically,
\begin{eqnarray*}
E_{-1} = {\frac {{a}^{2}\sqrt{\pi } \left( 256 \sqrt{2} -
51{a}^{2}b \right) }{2048 \sqrt{b^3}}}, \quad E_0 = \frac{3 a^2
\sqrt{2 \pi}}{8 \sqrt{b}},
\end{eqnarray*}
whereas the conserved quantities $E_1$ and $E_2$ are not expressed
in a closed form. Using numerical approximations of the integrals,
we determine the boundary of the well-posedness region in the
$(a,b)$-plane by finding the parameters $a$ and $b$ from the
condition $2\sqrt{2 E_1 E_2} = 1$. We also compute the boundary of
the wave breaking region in the $(a,b)$-plane by computing
$f_1$ and $f_2$ in (\ref{WB1}). Unlike the case of modulated pulses, we find
that the set $I_1$ is empty and $f_2$ may change the sign along the curve
on the $(a,b)$-plane. The two boundaries are shown on Figure \ref{fig:gauss},
where we can see that the two regions of global well-posedness and wave breaking are
disjoint.

Finally, we perform numerical simulations of the periodic Cauchy
problem (\ref{cauchy-periodic}) with the $1$-periodic initial data
\begin{equation}
u_0(x) = a \cos 2 \pi x, \quad a > 0.
\end{equation}
The two conserved quantities $E_0$ and $E_1$ in Lemma
\ref{lemma-Wayne} are computed analytically as
\begin{eqnarray*}
E_0 = \frac{1}{2} a^2, \quad E_1 = \frac{2}{\pi}
{\bf E}(2 \pi a i),
\end{eqnarray*}
where ${\bf E}$ stands for a complete elliptic integral.
Using the above conserved quantities we find
out that the sufficient condition for the wave breaking in Theorem
\ref{theorem-breaking-periodic} is satisfied for $a > 1.053$.

\begin{figure}
\begin{center}
\begin{tabular}{cc}
\includegraphics[width=0.4\textwidth]{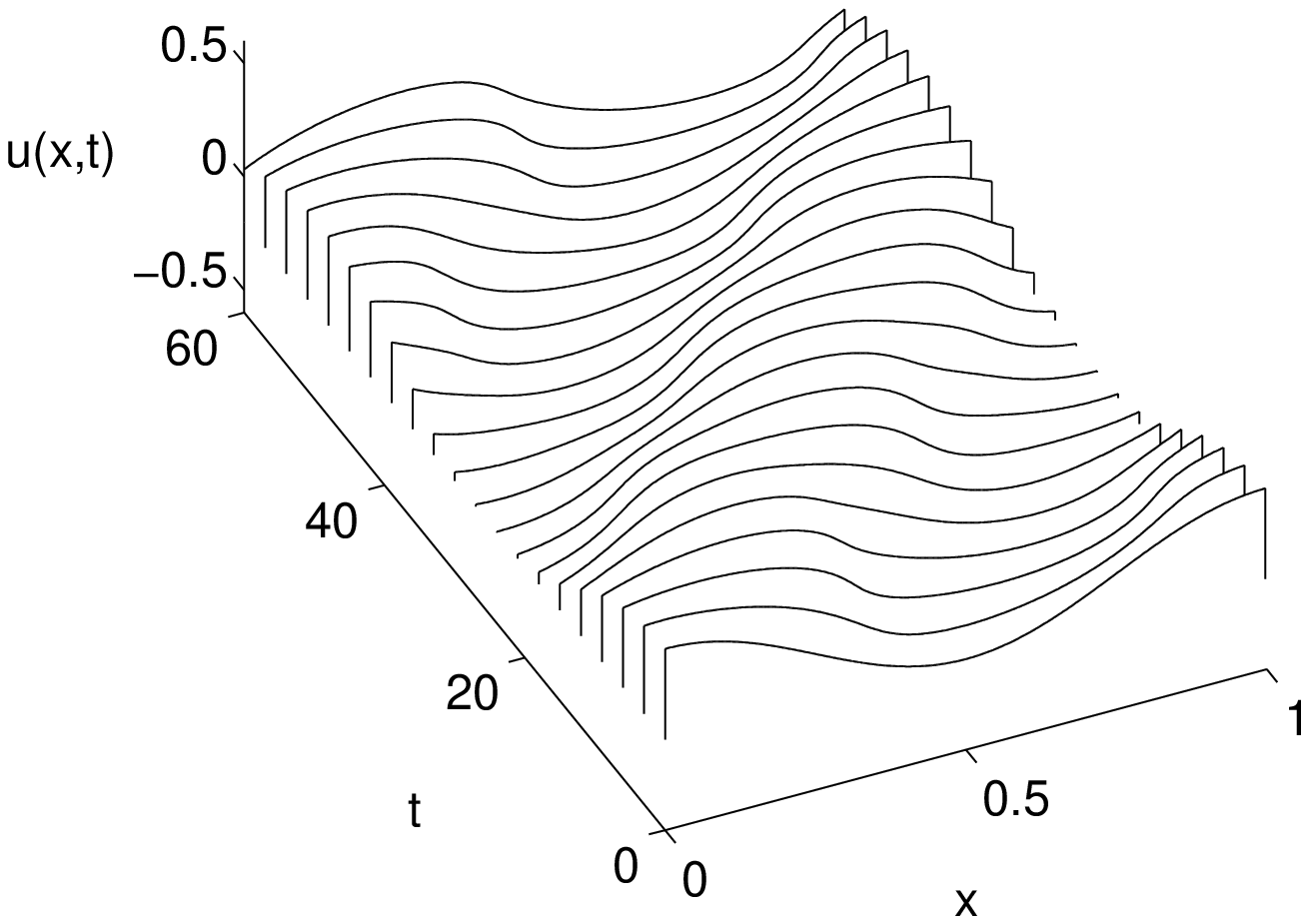} &
\includegraphics[width=0.4\textwidth]{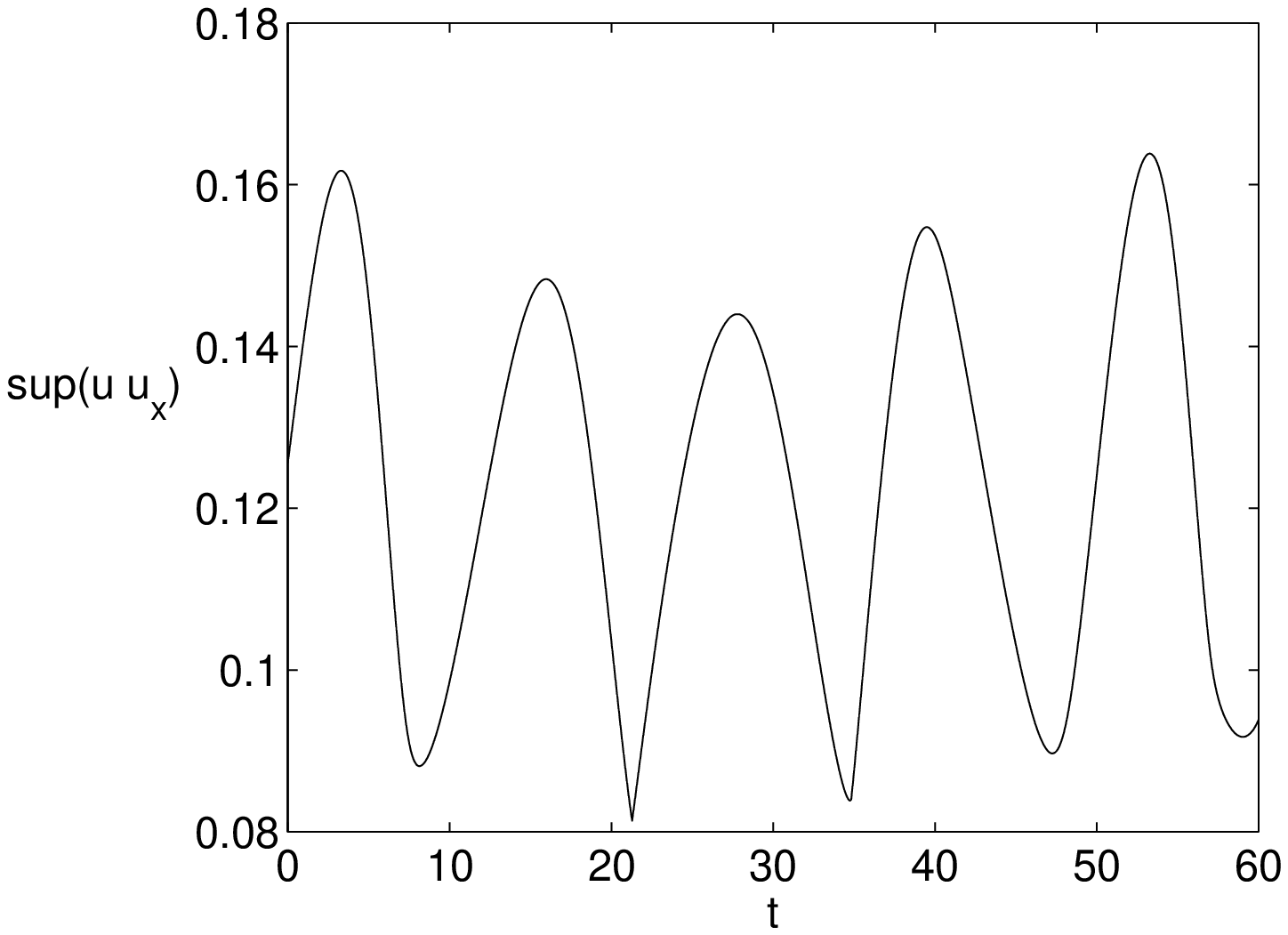} \\
\includegraphics[width=0.4\textwidth]{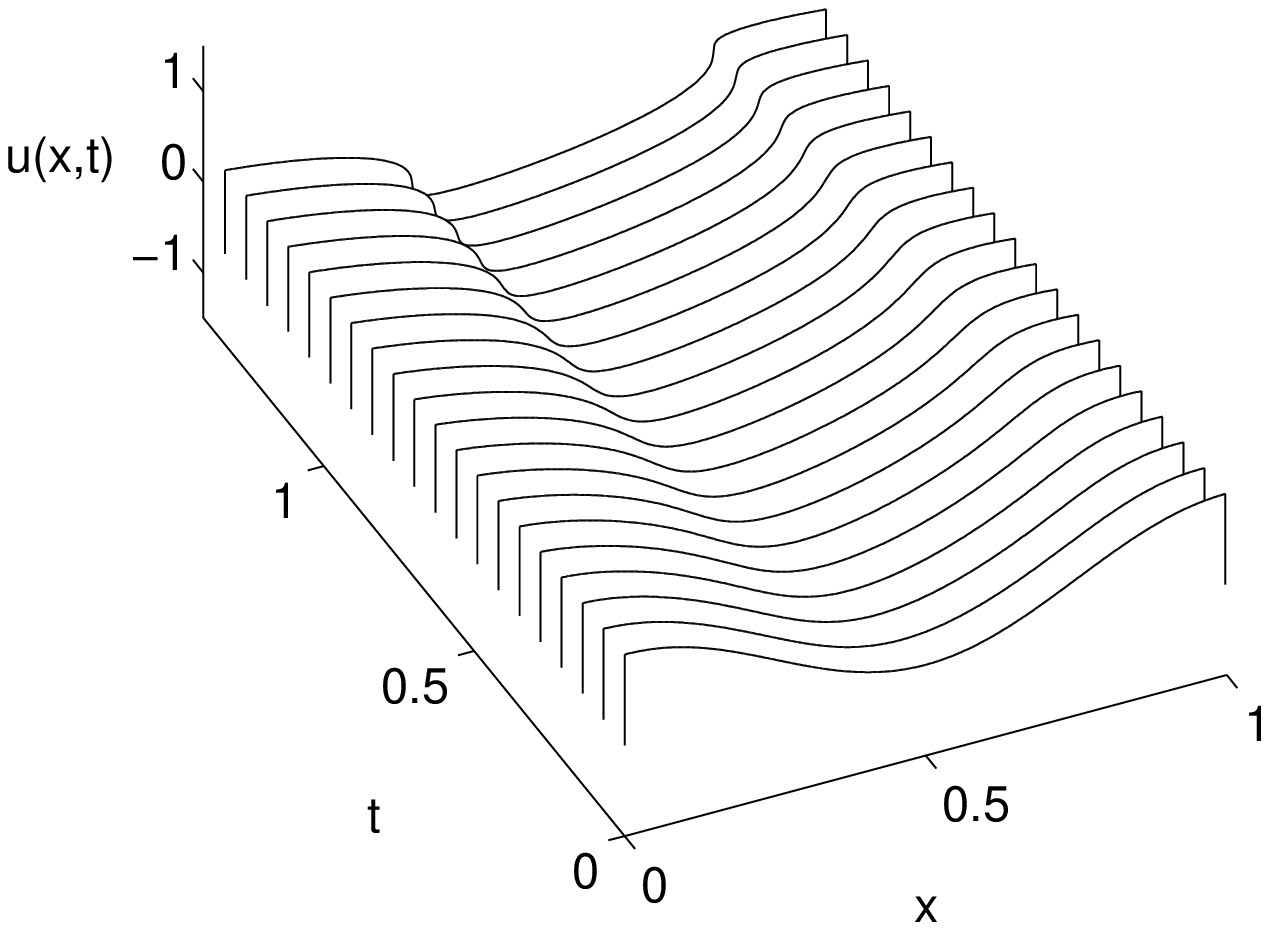} &
\includegraphics[width=0.4\textwidth]{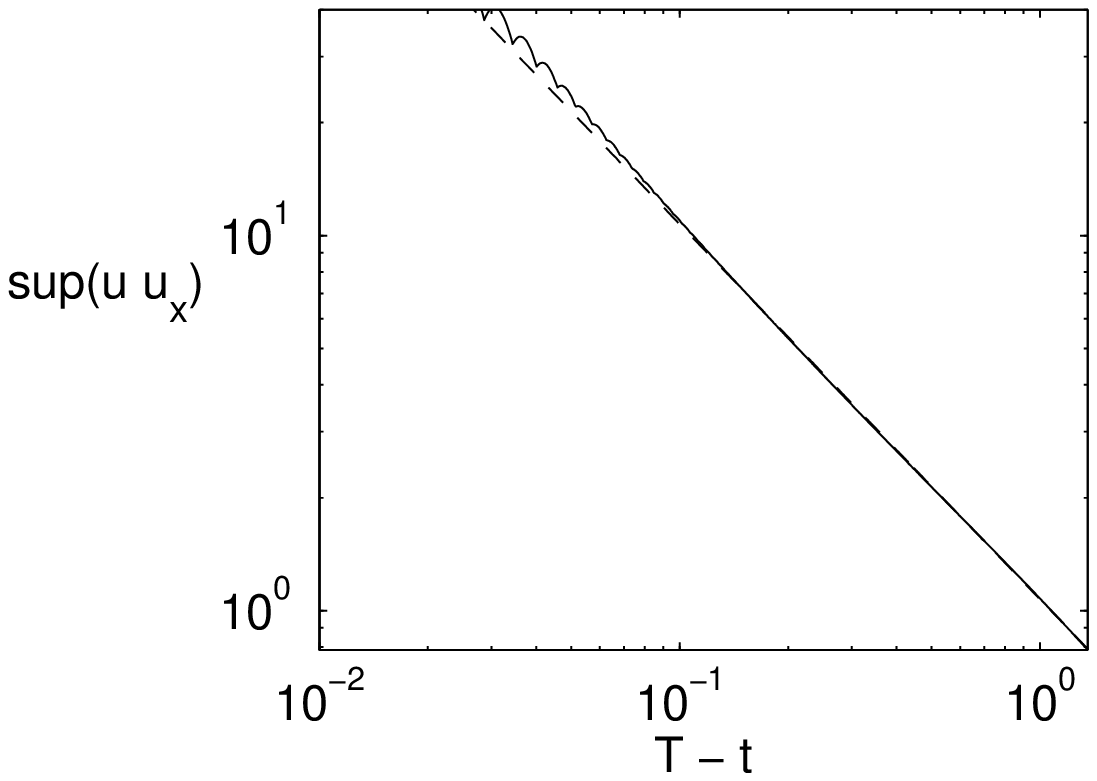}
\end{tabular}
\end{center}
\caption{Solution surface $u(x,t)$ (left) and the supremum norm
$W(t)$ (right) for $a=0.2$ (top) and $a=0.5$ (bottom). The dashed curve
on the bottom right picture shows the linear regression with $C=1.072$, $T=1.356$.}
\label{fig:periodic}
\end{figure}

\begin{figure}
\begin{center}
\includegraphics[width=0.4\textwidth]{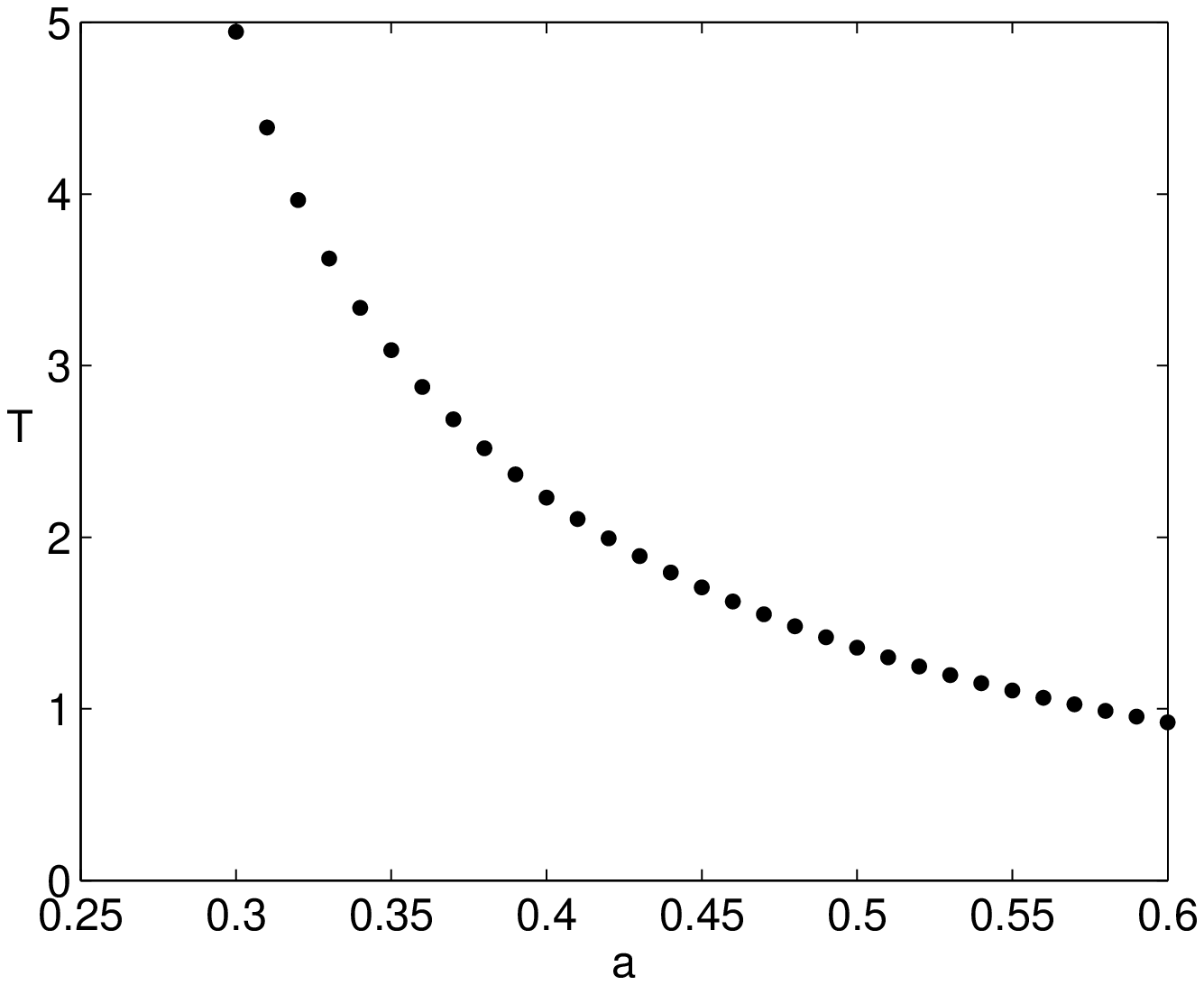}
\includegraphics[width=0.4\textwidth]{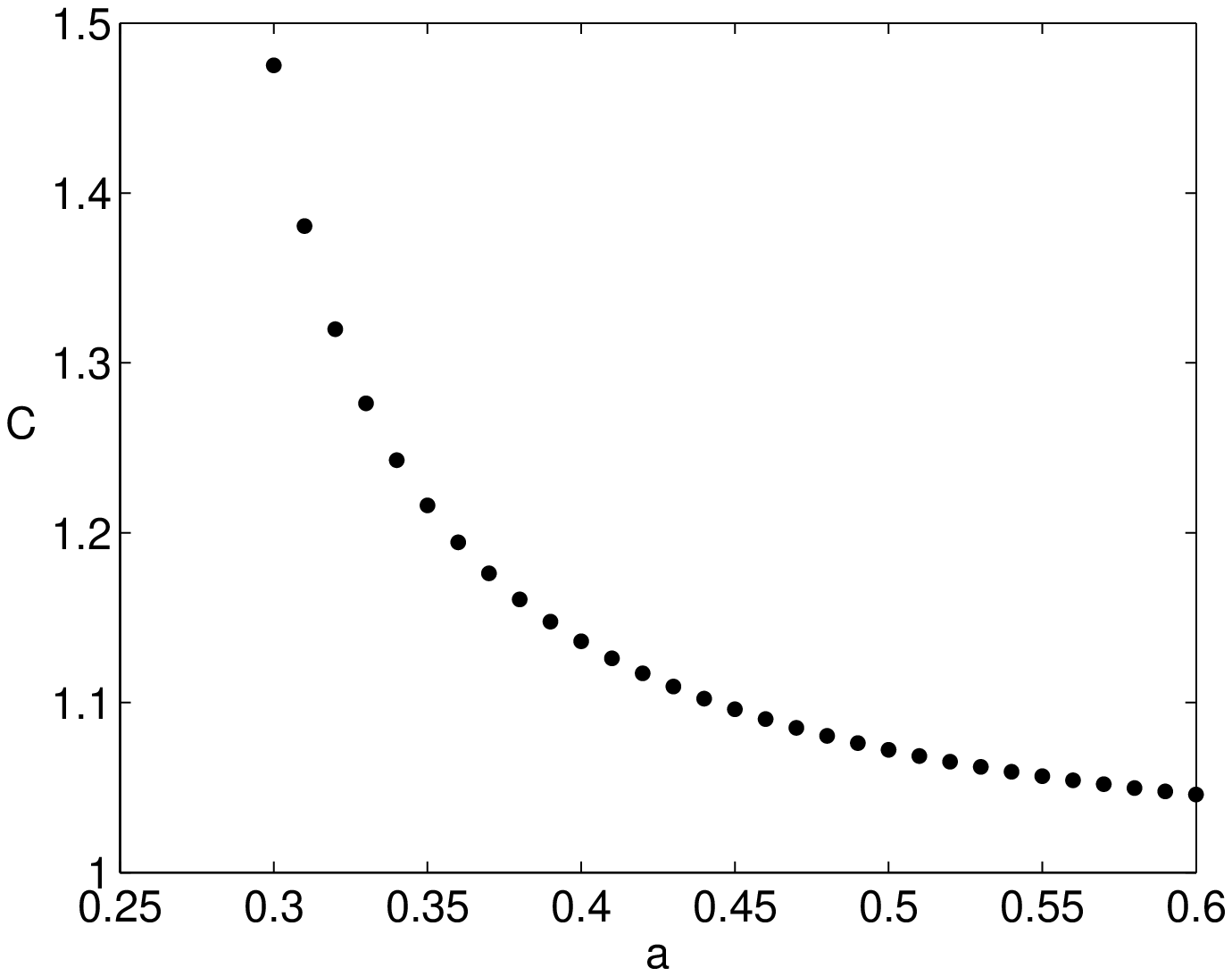}
\end{center}
\caption{Time of wave breaking $T$ versus $a$ (left). Constant $C$ of
the linear regression versus $a$ (right).}
\label{fig:T_vs_a}
\end{figure}

To illustrate the behaviour of a solution $u(x,t)$ to the Cauchy problem
(\ref{cauchy-periodic}), we perform numerical simulations using
a pseudospectral method. When the parameter $a$ is sufficiently
small, the value of
$$
W(t) := \sup_{x \in \mathbb{S}} u(x,t) u_x(x,t)
$$
remains bounded as shown on the top panel of Figure
\ref{fig:periodic} for $a=0.2$. On the other hand, when $a$
becomes larger, the wave breaking occurs, as on the bottom panel
of Figure \ref{fig:periodic} for $a = 0.5$. On the bottom right
panel of Figure \ref{fig:periodic} we show using the linear
regression that the curve $W^{-1}(t)$ is fitted well with 
the straight line $A + Bt$ for some coefficients $(A,B)$. 
Thus, we make a conclusion that
$$
W(t) \simeq \frac{C}{T - t} \quad \mbox{\rm for}
\quad 0 < T-t \ll 1,
$$
where $C = - B^{-1}$ and $T = -A B^{-1}$. Using the linear regression, we
also obtain pairs $(T,C)$ for different values of $a$. The results
are shown on Figure \ref{fig:T_vs_a}. Note that the constant $C$
approach $1$ as $a$ gets larger. This observation is consistent with the
exact blow-up law $W(t) = \frac{1}{T - t}$ obtained
for the dispersionless advection equation (\ref{simple-wave}) using the method of characteristics.

\end{document}